\documentclass[11pt,a4paper]{article}

\usepackage[utf8]{inputenc}
\usepackage[T1]{fontenc}
\usepackage{amsmath,amssymb,amsfonts}
\usepackage{amsthm}
\usepackage{geometry}
\usepackage{mathtools}
\usepackage{enumitem}
\usepackage{hyperref}
\usepackage{bm}
\usepackage{xcolor}

\newtheorem{lemma}{Lemma}

\newcommand{\bc}[1]{{\mathcal #1}}

\title{A note on the convergence of the eigenvalues in a subdomain to the continuous spectrum}
\author{Miroslav Bul\'{\i}\v{c}ek\thanks{Faculty of Mathematics and Physics, Charles University, 186 75 Prague 8, Czech Republic (email: mbul8060@karlin.mff.cuni.cz)} \and Bj\o rn Fredrik Nielsen\thanks{Faculty of Science and Technology, Norwegian University of Life Sciences, NO-1432 {\AA}s, Norway (email: bjorn.f.nielsen@nmbu.no).} \and Zden\v{e}k Strako\v{s}\thanks{Faculty of Mathematics and Physics, Charles University, 186 75 Prague 8, Czech Republic (email: strakos@karlin.mff.cuni.cz)}}
\begin{document}
\maketitle
\begin{abstract}
The SIGEST paper \cite{Nie24} characterized the spectrum of the preconditioned operator $\Delta^{-1}[\nabla \cdot (K\nabla u)]$ defined in a bounded and  open two dimensional domain $\Omega$, where $\Delta$ denotes the Laplacian and $K(x,y)$ is a continuous symmetric matrix-valued function. An important part of the analysis states that for a diagonal tensor $K$, which is constant in an open subdomain $S \subset \Omega $, the closed interval defined by its diagonal elements belongs to the spectrum of the preconditioned operator. This result is correct, but the proof that is present in \cite{Nie24} (and also in the original paper \cite{Ger20}) must be refined.

This paper presents the refined proof and extends the previous work in the following way. As shown in the cited papers, given any point $\lambda$ in the open interval defined by the elements of the diagonal tensor constant in $S$ and any point $(x_0, y_0) \in S$, a rectangular subdomain $\Sigma _l \subset S$ can be constructed such that the generalized eigenvalue problem associated with the preconditioned operator restricted to $\Sigma _l$ of arbitrarily small size has the eigenvalue $\lambda$ and infinitely many eigenfunctions. They are given by the solutions of {a locally defined wave equation}. However, this solution of the locally restricted  generalized eigenvalue problem cannot be extended to the whole domain $\Omega$. Using local rectangular subdomains with a size that {shrinks to zero}, the present paper instead constructs a Weyl singular sequence  of \emph{approximate} eigenfunctions associated with $\lambda$, which proves that $\lambda$ belongs to the spectrum of the preconditioned operator. Since self-adjoint operators in a separable Hilbert space can have at most {a countable set of eigenvalues}, this shows that the eigenvalues of the locally defined operator transform in the {limit to points of the continuous spectrum of the preconditioned operator defined throughout the entire domain}. 
\end{abstract}



\section{Introduction}
Using spectral decomposition, a symmetric 2D-tensor $K(x,y)$ can be written in the form 
\begin{equation*}
    K(x,y)=Q(x,y) \Lambda(x,y) Q^T(x,y), \quad (x,y) \in \Omega, 
\end{equation*}
where $Q$ is an orthogonal matrix and $\Lambda$ is a diagonal matrix. In \cite{Nie24} we analyze the spectrum of the preconditioned operator $\Delta^{-1} [\nabla \cdot (K\nabla u)]$ subject to homogeneous Dirichlet or Neumann boundary conditions, assuming that $K$ is continuous. Here,  $\Delta$ denotes the Laplacian, $\Omega \subset \mathbb{R}^2$ is a bounded and open Lipschitz domain and the function entries of $K$ are bounded and Lebesgue measurable. More specifically, we proved that the spectrum of the operator $\Delta^{-1} [\nabla \cdot (K\nabla u)]$, in a Sobolev space setting, is equal to the convex hull of the ranges of the diagonal function entries of $\Lambda$, discussed the consequences of this for operator preconditioning in the numerical treatment of elliptic PDEs and presented numerical experiments. 

An important step in the work \cite{Nie24} is the analysis of the case when $K$ is a diagonal tensor
\begin{equation} \label{eq:diag:tensor}
K(x,y)=\left[ 
\begin{array}{cc}
\kappa_1(x,y) & 0 \\
0 & \kappa_2(x,y)
\end{array}
\right],
\end{equation}
which is constant on an open subdomain $S \subset \Omega$. In particular, if the ranges of $\kappa_1(x,y)$ and $\kappa_2(x,y)$ over the domain $\Omega$ do not intersect, it is shown that the interval gap between them belongs to the spectrum of the preconditioned operator.

This issue is addressed by Lemma 3.2 in \cite{Nie24} and is the main concern of this paper: The statement of this lemma is correct, but its proof, as presented in \cite{Nie24}, suffers from some oversimplifications of how eigenvalues of locally defined operators  ``transform" to the spectrum of the associated globally defined operator. More precisely, simply constructing zero-extensions of locally defined eigenfunctions to the {entire domain $\Omega$ does not yield} a valid transition to the weak form of the generalized eigenvalue problem posed on all of $\Omega$, as certain line integrals do not vanish. 

The next section contains further details about these shortcomings. Section \ref{sec:main_result} outlines in subsection \ref{sec:proof} a correct proof of Lemma 3.2 in \cite{Nie24} and provides a detailed technical {analysis in subsection} \ref{sec:appendix}. The paper concludes with concluding remarks.

\section{Shortcomings} \label{sec:shortcomings}
Assume that the tensor $K(x,y)$ is diagonal and constant on the open subdomain $S \subset \Omega$,
\[
K(x,y) = \left[ 
\begin{array}{c c}
\overline{k}_1 & 0 \\
0 & \overline{k}_2
\end{array} 
\right],  \quad (x,y) \in S, 
\]
and consider the generalized eigenvalue problem 
\begin{equation}
\begin{split}
\label{eq:eigen}
\nabla \cdot (K \nabla u) &= \lambda \Delta u, \quad (x,y) \in \Omega, \\
u &= 0, \quad (x,y) \in\partial\Omega.
\end{split}
\end{equation}
{We will assume, without loss of generality, that $\overline{k}_1 < \overline{k}_2$. The case $\overline{k}_1 > \overline{k}_2$ can be handled analogously.} Denoting
\begin{equation}
\label{eq:constant_c}
c^2 = \frac{\lambda - \overline{k}_1}{\overline{k}_2 - \lambda} \,,
\end{equation}
{a solution $u$ of \eqref{eq:eigen}} must satisfy for any \(\lambda\) in the interval \((\overline{k}_1,\overline{k}_2)\), due to the indefiniteness of $K - \lambda I$ in $S$, 
\begin{equation} \label{eq:local-eigproblem}
    u_{yy}= c^2 \, u_{xx}, \quad (x,y) \in S, 
\end{equation}
i.e., $u$ must satisfy a wave equation on $S$. 

In connection with Lemma~3.2 in \cite{Nie24}, we therefore considered functions in the form 
\begin{equation*}
\label{B5.01}
\phi(x,y) =  \sin(n \pi  l^{-1} (x-x_0)) \sin(n \pi c l^{-1} (y-y_0)), 
\quad n \in \mathbb{N}, 
\end{equation*}
which solve the following Dirichlet problem for the wave equation 
\begin{equation}
\label{B5.1}
\begin{split}
\phi_{yy} &= c^2 \phi_{xx} \quad \mbox{in } \Sigma_l, \\
\phi & = 0 \quad \mbox{on } \partial \Sigma_l, 
\end{split}
\end{equation}
where $l$ is a positive constant that determines the size of the domain 
\[
\Sigma_l = (x_0,x_0 + l) \times (y_0,y_0+l/c). 
\]
Here, $(x_0,y_0)$ is a point in $S$, and $l$ is so small that $\Sigma_l \subset S$.  
It is important to note that $\Sigma_l$ depends on the choice of $c$ and therefore on $\lambda$.

We then erroneously claimed in \cite{Nie24} that the zero extension 
\[
v(x)
=
\begin{cases}
\phi(x), & x\in\Sigma_l,\\[4pt]
0, & x\notin\Sigma_l,
\end{cases}
\] 
is a weak solution of \eqref{eq:eigen}, i.e., that $v$ satisfies  
\begin{equation} \label{eq:weak_eigen}
    \int_\Omega (K \nabla v) \cdot \nabla \psi \, dx \, dy = \lambda \int_\Omega \nabla v \cdot \nabla \psi \,dx \, dy, \quad \forall \psi \in H^1_0(\Omega).  
\end{equation}

This is not correct: Since \( \phi \in H^1_0(\Sigma_l) \), the zero extension
\( v \) belongs to \(H^1_0(\Omega)\) and inside the subdomain \(\Sigma_l\), the function \( v \) 
satisfies the eigen-equation \eqref{eq:local-eigproblem}.
With a given fixed $\lambda$ and the restriction of (\ref{eq:eigen}) to $\Sigma_l$, the generalized eigenvalue $\lambda$ has infinite multiplicity.
However, for a general test function
\(\psi\in H^1_0(\Omega)\), an integration by parts over $\Omega$
yields
\[
\int_\Omega ( K  - \lambda I) \nabla v \cdot \nabla \psi \, dx \, dy
=
\int_{\partial\Sigma_l}
\bigl( (K-\lambda I)\nabla \phi \cdot \mathbf{n} \bigr)\,
\psi\,ds ,
\]
where \(\mathbf{n}\) denotes the outward unit normal on 
\(\partial\Sigma_l\).
This boundary term does not vanish in general, because 
\(\psi\) is arbitrary in \(H^1_0(\Omega)\) and is not required 
to vanish on \(\partial\Sigma_l\).
Consequently, the identity \eqref{eq:weak_eigen} does not hold. 

These considerations reveal that it is not possible to ``lift" local solutions of an eigen-equation on particular {subdomains constructed for the given 
$\lambda$} to global eigen-functions on the domain $\Omega$ using the zero extension. This is in accordance with the fact that the interval $(\overline{k}_1, \overline{k} _2)$ cannot be formed from the eigenvalues of the preconditioned {operator:} 
Indeed, the classical result states that a self-adjoint operator on a separable Hilbert space has at most a countable set of eigenvalues; see, e.g., \cite[Chapter~II, Section~8]{1955_JvN_Book}. 

Instead of focusing on localized solutions to the generalized eigenvalue problem (\ref{eq:eigen}), we will consider appropriate approximations of eigenfunctions leading to a singular Weyl sequence for each $\lambda \in (\overline{k}_1, \overline{k} _2)$.
The construction of this sequence takes into account the wave-like character (indefiniteness of $K - \lambda I$) of (\ref{eq:eigen}) in $S$, in the sense that the functions $\phi_l$ are chosen as oscillatory solutions reflecting the underlying dispersive structure of the operator. 
At the same time, these oscillatory profiles are progressively localized. More precisely, we start from a construction similar to the previous one, but allow the wave to evolve freely on a spatial scale that is asymptotically much larger than the localization region $\Sigma_l$, where  $l\to 0$. 
As a consequence, the resulting sequence of approximate eigenfunctions converges weakly to zero. However, due to concentration effects in the $H^1_0$-norm, the energy of the functions does not vanish. In particular, the $H^1_0$-norm remains bounded away from zero, while the norm of the residual of the generalized eigenvalue problem tends to zero in the limit.

\section{Main result} \label{sec:main_result}

%
Defining the operators $\mathcal{L}, \, \mathcal{A}: H_0^1(\Omega) \mapsto H^{-1}(\Omega)$ as
\begin{align}\label{eq:L}
& \langle \mathcal{L} \phi, \psi \rangle = \int_{\Omega} \nabla \phi  \cdot  \nabla \psi,  \quad \phi,  \psi \in H_0^1(\Omega), \\
\label{eq:A}
& \langle \mathcal{A} \phi, \psi \rangle = \int_{\Omega} K \nabla \phi  \cdot  \nabla \psi, \quad \phi,  \psi \in H_0^1(\Omega),
\end{align}
the paper \cite{Nie24} characterizes the spectrum 
\begin{equation}\label{eq:spectrum}
\mathrm{sp}(\mathcal{L}^{-1} \mathcal{A})
\equiv \left\{ \lambda \in \mathbb{C}; \,  \lambda \mathcal{I} - \mathcal{L}^{-1} \mathcal{A} \mbox{ does not have a bounded inverse} \right\}
\end{equation}
of the preconditioned operator
\begin{equation}\label{eq:operator}
\bc{L}^{-1}\bc{A}: H_0^1(\Omega)  \to H_0^1(\Omega).
\end{equation}
\noindent
As mentioned in the Introduction, this paper concerns the proof of Lemma 3.2 in \cite{Nie24}, which is restated as:

\begin{lemma} 
\label{lemma:constant}
Consider a diagonal tensor \eqref{eq:diag:tensor}
%
%
where the bounded and Lebesgue integrable functions $\kappa_i$, $i = 1,2$,  are constant on an open subdomain $S\subset\Omega$. Assuming that
\begin{equation}
\label{B6.1}
\sup_{(x,y)\in\Omega} \kappa_1(x,y) < \inf_{(x,y)\in\Omega}\kappa_2(x,y),
\end{equation}
the following closed interval belongs to the spectrum of $\bc{L}^{-1}\bc{A}$:
\begin{equation}\label{eq:constant}
\bigg[\sup_{(x,y)\in\Omega} \kappa_1(x,y), \, \inf_{(x,y)\in\Omega} \kappa_2(x,y)\bigg]  \subset \mathrm{sp}(\mathcal{L}^{-1} \mathcal{A}).
\end{equation}
The analogous statement obviously holds with the interchanging of the roles of $\kappa_1$ and $\kappa_2$ in (\ref{B6.1}) and (\ref{eq:constant}).
\end{lemma}

The construction in the proof presented in the following can be extended behind the statement of this lemma. However, we will prefer to remain in accordance with the logical structure of the previous work \cite{Nie24}.

\subsection{Outline of the proof} \label{sec:proof}
Consider an arbitrary fixed point $(x_0,y_0)\in S$. The tensor $K$ is in $S$ defined by the constants $\overline{k}_1$ and $\overline{k}_2$, 
\[
K(x,y) = \left[
\begin{array}{c c}
\overline{k}_1 & 0 \\
0 & \overline{k}_2
\end{array}
\right], \quad (x,y) \in S,
\]
where, by assumption,
\begin{align}\label{eq:con:proof}
\overline{k}_1 \leq \sup_{(x,y)\in\Omega} \kappa_1(x,y)< \inf_{(x,y)\in\Omega}\kappa_2(x,y) \leq \overline{k}_2.
\end{align} 

{The proof of Lemma \ref{lemma:constant} presented in this paper} constructs a Weyl singular sequence of approximate eigenfunctions. 
For an arbitrary fixed $\lambda$ in the interval $(\overline{k}_1,\overline{k}_2)$, we will construct a sequence of functions  $u_n(x,y) \in H_0^1(\Omega)$, {where $n$ will be a  sufficiently large natural number}, with the following properties:

\bigskip
\begin{description}
    \item[a)] The $H_0^1$-semi-norm is bounded from zero as $n \rightarrow \infty$, i.e., 
    \begin{equation*}
        | u_n |_{H^1(\Omega)} \geq \beta > 0, \quad \forall n. 
    \end{equation*}
    \item[b)] For functions $\psi\in B$, where $B:=\{\phi \in H^{1}_0(\Omega): \, |\phi|_{H^1(\Omega}\le 1\}$, 
    \begin{equation}\label{eq:con:Weyl}
        \lim_{n \rightarrow \infty} \,\sup_{\psi \in B} \,
        \langle \mathcal{A} u_n - \lambda \mathcal{L} u_n,\psi\rangle = 
        \lim_{n \rightarrow \infty} \, \sup_{\psi \in B} \, \int_{\Omega} (K - \lambda I)\,\nabla u_n \cdot \nabla \psi \, dx \, dy = 0.
    \end{equation}
    \noindent
    This can be equivalently written as
    \begin{equation*}
    \lim_{n \rightarrow \infty} \|(\mathcal{A}-\lambda\mathcal{L}){u}_n\|_{H^{-1}(\Omega)} = 0.
    \end{equation*}
\end{description}

\bigskip
\noindent
Assume now that {\bf a)} and {\bf b)} hold and that $\lambda$, as specified at the beginning of this paragraph, does not belong to the spectrum of $\bc{L}^{-1}\bc{A}$. Then $(\mathcal{L}^{-1} \mathcal{A}-\lambda \mathcal{I})^{-1}$ is bounded, which implies that 
\begin{equation*}
    (\mathcal{A}-\lambda\mathcal{L})^{-1} = (\mathcal{L}^{-1} \mathcal{A}-\lambda \mathcal{I})^{-1} \mathcal{L}^{-1}    
\end{equation*}
is bounded because $\mathcal{L}$ is continuously invertible. From this it follows that 
\begin{align*}
    | u_n |_{H^1(\Omega)} &= | (\mathcal{A}-\lambda\mathcal{L})^{-1} (\mathcal{A}-\lambda\mathcal{L}) u_n |_{H^1(\Omega)} \\
    &\leq \| (\mathcal{A}-\lambda\mathcal{L})^{-1} \| \, \|(\mathcal{A}-\lambda\mathcal{L}){u}_n\|_{H^{-1}(\Omega)} \\
    &\rightarrow 0
\end{align*}
as $n \rightarrow \infty$, which contradicts {\bf a)}. This shows that $\lambda$ indeed must be in the spectrum of $\bc{L}^{-1}\bc{A}$ and that  $(\overline{k}_1, \overline{k}_2) \subset \textnormal{sp}(\mathcal{L}^{-1}\mathcal{A})$; see also, e.g., \cite[Section~6.5]{1982Friedman_Book}. It remains to prove that if equality is achieved on any side of (\ref{eq:con:proof}), then the associated $\overline{k}_i$, $i=1$ and/or $i=2$ also belong to the spectrum of $\bc{L}^{-1}\bc{A}$. But this is trivially true because the spectrum is a closed set.

It remains to define the functions  $u_n (x, y)$ and prove that they obey {\bf a)} and {\bf b)}. 
Since the partial derivatives with respect to the variables $x$ respectively $y$ are in (\ref{eq:con:Weyl}) multiplied by the generally different values $\overline{k}_1$ respectively $\overline{k}_2$, we have to introduce a positive constant $c > 0$ that will deal with this geometrical asymmetry. As shown in the following section, its value will be given by (\ref{eq:constant_c}). {For sufficiently large values of $n$,} we consider the rectangle
\[
\Sigma_n := \left(x_0-\frac{1}{\sqrt{n}},\,x_0+\frac{1}{\sqrt{n}}\right)
\times
\left(y_0-\frac{1}{c\sqrt{n}},\,y_0+\frac{1}{c\sqrt{n}}\right)
\subset S
\]
that will reduce to the point $(x_0, y_0)$ as $n \rightarrow \infty$.
The function $u_n(x, y)$ determined below will vanish by construction in the complement $\Omega \backslash \Sigma_n$ and will belong to $H_0^1(\Omega)$. 

Denote the positive part of {a real number} by
\[
(z)_+ := \max\{z,0\},
\]
and for $n\in\mathbb{N}$ define
\[
f_n(t) := \sin(nt)\,(1 - n t^2)_+, \quad t\in\mathbb{R}
\]
that is equal to zero {outside the interval 
\[
\left( -\frac{1}{\sqrt{n}}, \frac{1}{\sqrt{n}} \right)
\]
and provides sufficient smoothness at its edges.} 
Finally, considering a sufficiently large $n$, for which $\Sigma_n \subset S$, 
define the function $u_n(x,y)$ as
\[
u_n(x,y) := \frac{1}{\sqrt{n}}\, f_n(x - x_0)\, f_n\big(c(y - y_0)\big),
\quad (x,y)\in \Omega.
\]
Observe that $u_n$ ``mainly" is a product of two sine modes. This construction is thus motivated by the discussion of the solutions of the local wave equations discussed in Section \ref{sec:shortcomings}.  

The proof of (\ref{eq:con:Weyl}) consists of three steps. First, the $H^1$ semi-norm
\[
|u_n|_{H^1(\Omega)}^2 =
\int_{\Omega} \big(|\partial_x u_n|^2 + |\partial_y u_n|^2\big)\,dx\,dy
\]
is expressed in the form
\[
|u_n|_{H^1(\Omega)}^2
=
\left(\frac{1}{c} + c\right)\frac{1}{n}
\left(\int_{-\frac{1}{\sqrt{n}}}^{\frac{1}{\sqrt{n}}}
|f_n(t)|^2\,dt\right)
\left(\int_{-\frac{1}{\sqrt{n}}}^{\frac{1}{\sqrt{n}}}
|f_n'(t)|^2\,dt\right).
\]
Using standard substitutions and the Riemann-Lebesgue lemma, this leads to 
the asymptotic result 
\[
\lim_{n \rightarrow \infty} |u_n|_{H^1(\Omega)}^2
= \left(\frac{1}{c} + c\right) \frac{1}{4} \int_{-1}^{1} (1 - s^{2})^{2}\, ds > 0
\]
that provides the desired {uniform bound {\bf a)}} from zero for all large $n$.

Second, the integral 
\begin{equation}\label{eq:con:Weyl:In}
I_n = \int_{\Sigma_n} (K - \lambda I)\,\nabla u_n \cdot \nabla \psi \, dx \, dy    
\end{equation}
in (\ref{eq:con:Weyl}) will be rewritten as 
\begin{align*}
I_n
&=
\frac{\overline{k}_1 - \lambda}{\sqrt{n}}
\int_{\Sigma_n}
f_n'(x - x_0)\, f_n\!\bigl(c(y - y_0)\bigr)\,
\partial_x \psi(x,y)\, dx\, dy
\\
&\quad
+
\frac{c(\overline{k}_2 - \lambda)}{\sqrt{n}}
\int_{\Sigma_n}
f_n(x - x_0)\, f_n'\!\bigl(c(y - y_0)\bigr)\,
\partial_y \psi(x,y)\, dx\, dy .
\end{align*}
After substituting for $f_n$ it will get the form (we will skip the {arguments $x, \, y$} in the function $\psi(x,y)$)
\begin{align*}
I_n =
\frac{\overline{k}_1 - \lambda}{\sqrt{n}}
&
\int_{\Sigma_n}
\Bigl[
n \cos\bigl(n(x - x_0)\bigr)\bigl(1 - n(x - x_0)^2\bigr)
- 2n(x - x_0)\sin\bigl(n(x - x_0)\bigr)
\Bigr] \\
 & \qquad \qquad \qquad \cdot\, f_n\!\bigl(c(y - y_0)\bigr)\,
\partial_x \psi\, dx\, dy
\\
+
\frac{c(\overline{k}_2 - \lambda)}{\sqrt{n}}
&
\int_{\Sigma_n}
\Bigl[
n \cos\bigl(nc(y - y_0)\bigr)\bigl(1 - n c^2 (y - y_0)^2\bigr)
- 2nc(y - y_0)\sin\bigl(nc(y - y_0)\bigr)
\Bigr] \\
 & \qquad \qquad \qquad \cdot \, f_n(x - x_0)\,
\, \partial_y \psi \, dx\, dy \\
&=: I_n^{\sin} + I_n^{\cos},
\end{align*}
where 
\begin{align*}
I_n^{\sin}
&=
- \frac{2n}{\sqrt{n}}
\int_{\Sigma_n}
(\overline{k}_1 - \lambda)(x - x_0)\sin\bigl(n(x - x_0)\bigr)
 \cdot\,
f_n\!\bigl(c(y - y_0)\bigr)\,\partial_x \psi \, dx\, dy 
\\
&\quad
- \frac{2n c^2}{\sqrt{n}}
\int_{\Sigma_n}
(\overline{k}_2 - \lambda)(y - y_0)\sin\bigl(nc(y - y_0)\bigr)
 \cdot\,
f_n(x - x_0)\,\partial_y \psi\, dx\, dy 
\end{align*}
and
\begin{align*}
I_n^{\cos}
&=
\frac{n}{\sqrt{n}}
\int_{\Sigma_n}
(\overline{k}_1 - \lambda)
\cos\bigl(n(x - x_0)\bigr)
\bigl(1 - n(x - x_0)^2\bigr)
 \cdot\,
f_n\!\bigl(c(y - y_0)\bigr)\,\partial_x \psi\, dx \, dy
\\
&+
\frac{nc}{\sqrt{n}}
\int_{\Sigma_n}
(\overline{k}_2 - \lambda)
\cos\bigl(nc(y - y_0)\bigr)
\bigl(1 - nc^2 (y - y_0)^2\bigr)
 \cdot\,
f_n(x - x_0)\,\partial_y \psi \, dx \,dy .
\end{align*}

\noindent

Third, it is proved that both 
$$\sup_{| \psi |_{H^1(\Omega)} = 1} I_n^{\sin}  \rightarrow 0 \quad \mbox{and} \quad \sup_{| \psi |_{H^1(\Omega)} = 1} I_n^{\cos} \rightarrow 0$$ 
as $n \rightarrow \infty$, {which completes the argument for {\bf b)}}. 
The technical steps require setting the constant $c$ as (see (\ref{eq:constant_c}))
\[
c^2 = \frac{\lambda - \overline{k}_1}{\overline{k}_2 - \lambda}. 
\]
Detailed derivation is given  below.

\subsection{Technical details of the derivation} \label{sec:appendix}

We will proceed according to the three steps outlined above.

\paragraph{{\bf $\mathbf{H_0^1}$ semi-norm of $\mathbf{u_n(x, y)}$}} We are going to evaluate
\[
|u_n|_{H^1(\Omega)}^2
=
\int_{\Sigma_n} \big(|\partial_x u_n|^2 + |\partial_y u_n|^2\big)\,dx\,dy.
\]
Considering the $x$-derivative
\[
\partial_x u_n(x,y)
=
\frac{1}{\sqrt{n}}\, f_n'(x-x_0)\, f_n\big(c(y-y_0)\big),
\]
we get
\[
\int_{\Sigma_n} |\partial_x u_n|^2\,dx\,dy
=
\frac{1}{n}
\int_{x_0-\frac{1}{\sqrt{n}}}^{x_0+\frac{1}{\sqrt{n}}}
|f_n'(x-x_0)|^2\,dx
\int_{y_0-\frac{1}{c\sqrt{n}}}^{y_0+\frac{1}{c\sqrt{n}}}
|f_n(c(y-y_0))|^2\,dy,
\]
which with the substitutions $t=x-x_0$ and $s=c(y-y_0)$ results in
%
%
\[
\int_{\Sigma_n} |\partial_x u_n|^2\,dx\,dy
=
\frac{1}{c\,n}
\left(\int_{-\frac{1}{\sqrt{n}}}^{\frac{1}{\sqrt{n}}}
|f_n'(t)|^2\,dt\right)
\left(\int_{-\frac{1}{\sqrt{n}}}^{\frac{1}{\sqrt{n}}}
|f_n(s)|^2\,ds\right).
\]
Analogously, for the $y$-derivative
\[
\partial_y u_n(x,y)
=
\frac{c}{\sqrt{n}}\, f_n(x-x_0)\, f_n'\big(c(y-y_0)\big)
\]
we get
\[
\int_{\Sigma_n} |\partial_y u_n|^2\,dx\,dy
=
\frac{c^2}{n}
\int_{x_0-\frac{1}{\sqrt{n}}}^{x_0+\frac{1}{\sqrt{n}}}
|f_n(x-x_0)|^2\,dx
\int_{y_0-\frac{1}{c\sqrt{n}}}^{y_0+\frac{1}{c\sqrt{n}}}
|f_n'(c(y-y_0))|^2\,dy.
\]
Using the same substitutions as above and noticing {that}
%
%
\[
\int_{y_0-\frac{1}{c\sqrt{n}}}^{y_0+\frac{1}{c\sqrt{n}}}
|f_n'(c(y-y_0))|^2\,dy
=
\frac{1}{c}
\int_{-\frac{1}{\sqrt{n}}}^{\frac{1}{\sqrt{n}}}
|f_n'(s)|^2\,ds
\]
results in
\[
\int_{\Sigma_n} |\partial_y u_n|^2\,dx\,dy
=
\frac{c}{n}
\left(\int_{-\frac{1}{\sqrt{n}}}^{\frac{1}{\sqrt{n}}}
|f_n(t)|^2\,dt\right)
\left(\int_{-\frac{1}{\sqrt{n}}}^{\frac{1}{\sqrt{n}}}
|f_n'(s)|^2\,ds\right).
\]
Combining both contributions and simplifying the notation,
\begin{equation}\label{eq:midterm}
|u_n|_{H^1(\Omega)}^2
=
\left(\frac{1}{c} + c\right)\frac{1}{n}
\left(\int_{-\frac{1}{\sqrt{n}}}^{\frac{1}{\sqrt{n}}}
|f_n(t)|^2\,dt\right)
\left(\int_{-\frac{1}{\sqrt{n}}}^{\frac{1}{\sqrt{n}}}
|f_n'(t)|^2\,dt\right).
\end{equation}

Additional change of variables 
\[ t=\frac{s}{\sqrt{n}}, \quad \mbox{i.e.,} \quad nt=\sqrt{n}\,s,\quad nt^2=s^2,\quad s\in[-1,1],\]
and substitution for the function $f(t)$, provides
\[
\begin{aligned}
&\int_{-\frac{1}{\sqrt{n}}}^{\frac{1}{\sqrt{n}}} |f_n(t)|^2\,dt
=
\frac{1}{\sqrt{n}}
\int_{-1}^{1}
\sin^2(\sqrt{n}\,s)(1-s^2)^2\,ds,
\\
&\int_{-\frac{1}{\sqrt{n}}}^{\frac{1}{\sqrt{n}}} |f_n'(t)|^2\,dt
=
\frac{1}{\sqrt{n}}
\int_{-1}^{1}
\left[
n\cos(\sqrt{n}\,s)(1-s^2)
-
2\sqrt{n}\,s\,\sin(\sqrt{n}\,s)
\right]^2 ds \\
 &\hspace{2.75cm}=\frac{n^2}{\sqrt{n}}
\int_{-1}^{1}
\left[
\cos(\sqrt{n}\,s)(1-s^2)
-
\frac{2}{\sqrt{n}} \,s\, \sin(\sqrt{n}\,s)
\right]^2 ds.
\\[0.5em]
\end{aligned}
\]
The powers of $n$ will cancel out {in the multiplication of} the terms in (\ref{eq:midterm})). Using the goniometric formula and the Riemann–Lebesgue lemma,
\begin{align*}
\lim_{n \rightarrow \infty} 
&\int_{-1}^{1} \sin^{2}(\sqrt{n}\, s)\,(1 - s^{2})^{2}\, ds \\
&= \frac{1}{2}  \int_{-1}^{1} (1 - s^{2})^{2}\, ds 
      - \frac{1}{2}  \, \lim_{n \rightarrow \infty} \int_{-1}^{1} (1 - s^{2})^{2} \cos(2\sqrt{n}\, s)\, ds 
 = \frac{1}{2} \int_{-1}^{1} (1 - s^{2})^{2}\, ds .
\end{align*}
Concerning the term with the derivative, 
\begin{align*}
    \lim_{n \rightarrow \infty}  
    &\int_{-1}^{1}\left[ 
    \cos(\sqrt{n}\,s)(1-s^2) -
    \frac{2}{\sqrt{n}}\,s\,\sin(\sqrt{n}\,s)
    \right]^2 ds \\
    &= 
    \lim_{n \rightarrow \infty}  
    \int_{-1}^{1} \cos^2(\sqrt{n}\,s)(1-s^2)^2 \, ds \\
    &=  \frac{1}{2} \int_{-1}^{1} (1 - s^{2})^{2}\, ds
+ \frac{1}{2} \, \lim_{n \rightarrow \infty} \int_{-1}^{1} (1 - s^{2})^{2} \cos(2\sqrt{n}\, s)\, ds 
     = \frac{1}{2} \int_{-1}^{1} (1 - s^{2})^{2}\, ds.
\end{align*} 
We conclude that
\[
\lim_{n \rightarrow \infty} |u_n|_{H^1(S)}^2
= \left(\frac{1}{c} + c\right) \frac{1}{4} \int_{-1}^{1} (1 - s^{2})^{2}\, ds > 0,
\]
which finishes the first part of the proof, {i.e., the proof of {\bf a)}}.

\bigskip
\paragraph{{\bf Convergence $\mathbf{\sup_{\psi \in B} \, I_n^{\sin} \rightarrow 0}$}}

Substitutions
\[
t = x-x_0, s = c(y-y_0), \widetilde{\psi}(t,s) = \psi (x_0 + t, y_0 + \frac{1}{c} s),
\partial_y \psi = c \, \partial_s \widetilde{\psi}(t,s),
\]
will transfer the integration domain to 
\[
\widetilde{\Sigma}_n
=
\left(
- \frac{1}{\sqrt{n}},\;
\frac{1}{\sqrt{n}}
\right)
\times
\left(
- \frac{1}{\sqrt{n}},\;
\frac{1}{\sqrt{n}}
\right)
\]
and provide
\begin{align*}
I_n^{\sin}
&=
- \frac{2\sqrt{n}}{c}
\int_{-\frac{1}{\sqrt{n}}}^{\frac{1}{\sqrt{n}}} \int_{-\frac{1}{\sqrt{n}}}^{\frac{1}{\sqrt{n}}}
(\overline{k}_1 - \lambda) \, t \sin(nt)
 \cdot\,
f_n(s)\,\partial_t \widetilde{\psi}(t,s)\, dt\,ds
\\
&\quad
- 2\sqrt{n} \, c
\int_{-\frac{1}{\sqrt{n}}}^{\frac{1}{\sqrt{n}}} \int_{-\frac{1}{\sqrt{n}}}^{\frac{1}{\sqrt{n}}}\,
(\overline{k}_2 - \lambda) \, s \sin(ns)
 \cdot\,
f_n(t)\,\partial_s \widetilde{\psi}(t,s)\, dt\, ds .
\end{align*}
Using $c$ as defined in (\ref{eq:constant_c}),
\[
c^2 = \frac{\lambda - \overline{k}_1}{\overline{k}_2 - \lambda} ,
\]
and the notation
\begin{equation*}
    d =\sqrt{(\lambda - \overline{k}_1)(\overline{k}_2 - \lambda)} \,,
\end{equation*}
we find that
\begin{align*}
I_n^{\sin}
&=
2\sqrt{n} \, d
\int_{-\frac{1}{\sqrt{n}}}^{\frac{1}{\sqrt{n}}} \int_{-\frac{1}{\sqrt{n}}}^{\frac{1}{\sqrt{n}}}
t\, \sin (nt)
 \cdot\,
f_n (s)\,\partial_t \widetilde{\psi}(t,s)\, dt\, ds
\\
&\quad  \quad
- 2\sqrt{n} \, d
\int_{-\frac{1}{\sqrt{n}}}^{\frac{1}{\sqrt{n}}} \int_{-\frac{1}{\sqrt{n}}}^{\frac{1}{\sqrt{n}}}
s\,\sin(n s)
 \cdot\,
f_n(t)\,\partial_s \widetilde{\psi}(t,s)\, dt\, ds \\
&= 2\sqrt{n} \, d 
\int_{-\frac{1}{\sqrt{n}}}^{\frac{1}{\sqrt{n}}} \int_{-\frac{1}{\sqrt{n}}}^{\frac{1}{\sqrt{n}}}
\Bigl(g(t,s),-g(s,t)\Bigr) \cdot \nabla \widetilde{\psi}(t,s) \, dt\, ds,
\end{align*}
where
 $$g(t,s) = t\, \sin (nt) \cdot\,f_n (s).$$
Consequently,
\begin{align*}
  (I_n^{\sin})^2 &\leq 8n {d^2} \,| \widetilde{\psi} |_{H^1(\Omega)}^2 
  \int_{- \frac{1}{\sqrt{n}}}^{\frac{1}{\sqrt{n}}}
\int_{- \frac{1}{\sqrt{n}}}^{\frac{1}{\sqrt{n}}} t^2 \sin^2(nt) f_n^2(s) \, dt \, ds \\
&=8 \sqrt{n} {d^2} \, | \widetilde{\psi} |_{H^1(\Omega)}^2 
  \int_{- \frac{1}{\sqrt{n}}}^{\frac{1}{\sqrt{n}}}
 (t \, \sin(nt))^2 \, dt \,  \int_{- \frac{1}{\sqrt{n}}}^{\frac{1}{\sqrt{n}}} \sqrt{n} f_n^2(s) \, ds
\end{align*}
Substituting $\sqrt{n}t=x$ and $\sqrt{n}s=x$, respectively, yields 
\begin{align*}
    (I_n^{\sin})^2 &\leq 8 \sqrt{n} d^2 \, | \widetilde{\psi} |_{H^1(\Omega)}^2 
  \int_{-1}^{1}
 \frac{1}{n \sqrt{n}}(x \sin(\sqrt{n}x))^2 \, dx \,  
 \int_{-1}^{1} \sqrt{n} f_n^2(x/\sqrt{n}) \, \frac{1}{\sqrt{n}} \, dx \\
 &= \frac{8 d^2} {n} \, | \widetilde{\psi} |_{H^1(\widetilde{\psi})}^2
  \int_{-1}^{1}
 (x \sin(\sqrt{n}x))^2 \, dx \,  \int_{-1}^{1} f_n^2(x/\sqrt{n}) \,  dx  .\\
\end{align*}
We conclude that 
\begin{equation*}
\lim_{n \rightarrow \infty}  \sup_{| \psi |_{H^1(\Omega)} = 1} I_n^{\sin} = 0. 
\end{equation*}
Here we have used $|{\psi} |_{H_0^1(\Omega)} = |\widetilde{\psi} |_{H_0^1(\Omega)} / c $.

\bigskip
\paragraph{{\bf Convergence $\mathbf{\sup_{\psi \in B} \, I_n^{\cos} \rightarrow 0}$}} Using, as above,
\[
t = x-x_0, s = c(y-y_0), \widetilde{\psi}(t,s) = \psi (x_0 + t, y_0 + \frac{1}{c} s),
\partial_y \psi = c \, \partial_s \widetilde{\psi}(t,s),
\]
\begin{align*}
I_n^{\cos}
&=
\frac{\sqrt{n}}{c}
\int_{-\frac{1}{\sqrt{n}}}^{\frac{1}{\sqrt{n}}}
\int_{-\frac{1}{\sqrt{n}}}^{\frac{1}{\sqrt{n}}}
\,
(\overline{k}_1 - \lambda)
\cos(nt)(1 - nt^2) \cdot\,
f_n(s)\,\partial_t \widetilde{\psi}(t,s)\, dt\, ds
\\
&\quad
+
\sqrt{n} \, c
\int_{-\frac{1}{\sqrt{n}}}^{\frac{1}{\sqrt{n}}}
\int_{-\frac{1}{\sqrt{n}}}^{\frac{1}{\sqrt{n}}}
\,
(\overline{k}_2 - \lambda)
\cos(ns)(1 - ns^2)
\cdot\,
f_n(t)\,\partial_s \widetilde{\psi}(t,s)\, dt\, ds, 
\end{align*}
which can be expressed in the form 
\begin{align*}
I_n^{\cos}
&=
- \sqrt{n} \,d
\int_{-\frac{1}{\sqrt{n}}}^{\frac{1}{\sqrt{n}}}
\int_{-\frac{1}{\sqrt{n}}}^{\frac{1}{\sqrt{n}}}
\,
\cos(nt)\,(1 - nt^2) \,\cdot\,
\sin(ns)\,(1 - ns^2)\,
\partial_t \widetilde{\psi}(t, s)\, dt\, ds
\\
&\quad
+
\sqrt{n} \,d
\int_{-\frac{1}{\sqrt{n}}}^{\frac{1}{\sqrt{n}}}
\int_{-\frac{1}{\sqrt{n}}}^{\frac{1}{\sqrt{n}}}
\,
\cos(ns)\,(1 - ns^2) \,\cdot\,
\sin(nt)\,(1 - nt^2)\,
\partial_s \widetilde{\psi}(t, s)\, dt\, ds .
\end{align*}
We will express{, for later convenience,} the sine functions as the derivative of cosines
and slightly rearrange the integrands,
\begin{align*}
I_n^{\cos}
&=
\frac{d}{\sqrt{n}}
\int_{-\frac{1}{\sqrt{n}}}^{\frac{1}{\sqrt{n}}}
\int_{-\frac{1}{\sqrt{n}}}^{\frac{1}{\sqrt{n}}}
\Bigl[
\Bigl(\cos(nt)\,(1 - nt^2)\,(1 - ns^2)\,\partial_t \widetilde{\psi}\Bigr)\,
\partial_s \cos(ns)\bigr)\,
\\
&\qquad\qquad
-
\Big(\cos(ns)\,(1 - ns^2)\,(1 - nt^2)\,\partial_s \widetilde{\psi} \Bigr) \,
\partial_t \cos(nt)
\Bigr]
\, dt\, ds .
\end{align*}
Now we use integration by parts in both variables to move the
derivatives off the cosine factors. We start with the first part of 
the integrand that will change with integration  by parts in the variable $s$ as
\begin{align*}
\Bigl(\cos(nt)\,(1 - nt^2)\,(1 - ns^2)\,\partial_t \widetilde{\psi}\Bigr)\, 
&
\partial_s \cos(ns) \rightarrow \\
&- \cos(ns) \, \partial_s \Bigl(\cos(nt)\,(1 - nt^2)\,(1 - ns^2)\,\partial_t \widetilde{\psi}\Bigr), 
\end{align*}
where the boundary terms vanish because
\[
(1 - ns^2) = 0
\quad \text{at } \quad s = \pm \frac{1}{\sqrt{n}}.
\]
Hence, the first term becomes
\[
- \cos(nt)\,\cos(ns)\,
\partial_s \Bigl( (1 - nt^2)(1 - ns^2)\,\partial_t \widetilde{\psi} \Bigr).
\]
Analogously for the second term with integration by parts in the variable $t$,
\begin{align*}
- \Bigl(\cos(ns)\,(1 - ns^2)\,(1 - nt^2)\,\partial_s \widetilde{\psi}\Bigr)\, 
&
\partial_t \cos(nt) \rightarrow \\
&
\cos(nt) \, \partial_t \Bigl(\cos(ns)\,(1 - ns^2)\,(1 - nt^2)\,\partial_s \widetilde{\psi}\Bigr), 
\end{align*}
where again the boundary terms vanish because
\[
(1 - nt^2) = 0
\quad \text{at } \quad t = \pm \frac{1}{\sqrt{n}}.
\]
Thus the second term becomes
\[
+ \cos(nt)\,\cos(ns)\,
\partial_t \Bigl( (1 - ns^2)(1 - nt^2)\,\partial_s \widetilde{\psi} \Bigr).
\]
We will therefore get{, after the integration by parts,}
\begin{align*}
I_n^{\cos}
&=
\frac{d}{\sqrt{n}}
\int_{-\frac{1}{\sqrt{n}}}^{\frac{1}{\sqrt{n}}}
\int_{-\frac{1}{\sqrt{n}}}^{\frac{1}{\sqrt{n}}}
\cos(nt)\,\cos(ns)
\\
&\qquad \cdot
\Bigl[
- \partial_s \Bigl( (1 - nt^2)(1 - ns^2)\,\partial_t \widetilde{\psi} \Bigr)
+ \partial_t \Bigl( (1 - ns^2)(1 - nt^2)\,\partial_s \widetilde{\psi} \Bigr)
\Bigr]
\, dt\, ds ,
\end{align*}
and, using the fact that 
$\partial_s \partial_t \widetilde{\psi} = \partial_t \partial_s \widetilde{\psi}$,
$$
I_n^{\cos} =
\frac{d}{\sqrt{n}}
\int_{-\frac{1}{\sqrt{n}}}^{\frac{1}{\sqrt{n}}}
\int_{-\frac{1}{\sqrt{n}}}^{\frac{1}{\sqrt{n}}}
\cos(nt)\,\cos(ns)
\cdot
\Bigl[
2ns\,(1 - nt^2)\,\partial_t \widetilde{\psi}
-
2nt\,(1 - ns^2)\,\partial_s \widetilde{\psi}
\Bigr]
\, dt\, ds .
$$
This finally gives, 
\begin{align*}
(I_n^{\cos})^2
&\leq
\frac{2d^2}{n} \,
|\widetilde{\psi}|_{H^1(\Omega)}^2
\int_{-\frac{1}{\sqrt{n}}}^{\frac{1}{\sqrt{n}}}
\int_{-\frac{1}{\sqrt{n}}}^{\frac{1}{\sqrt{n}}}
(\cos(nt))^2\,(\cos(ns))^2 \, 4(ns)^2 \,(1 - nt^2)^2 \, dt\, ds \\
&= \frac{8d^2}{n}\,
|\widetilde{\psi}|_{H^1(\Omega)}^2 
\int_{-\frac{1}{\sqrt{n}}}^{\frac{1}{\sqrt{n}}} (\cos(nt))^2\,(1 - nt^2)^2 \, dt \, 
\int_{-\frac{1}{\sqrt{n}}}^{\frac{1}{\sqrt{n}}} (\cos(ns))^2 \, (ns)^2 \, ds \\
&\leq \frac{8d^2}{n}\,
|\widetilde{\psi}|_{H^1(\Omega)}^2 
\int_{-\frac{1}{\sqrt{n}}}^{\frac{1}{\sqrt{n}}} (1 - nt^2)^2 \, dt \, 
\int_{-\frac{1}{\sqrt{n}}}^{\frac{1}{\sqrt{n}}} (ns)^2 \, ds \\
&= \frac{8d^2}{n}\,
|\widetilde{\psi}|_{H^1(\Omega)}^2 \frac{16}{15\sqrt{n}} \, \frac{2}{3}\sqrt{n}, 
\end{align*}
{and we can conclude that}  
\begin{equation*}
\lim_{n \rightarrow \infty}  \sup_{| \psi |_{H^1(\Omega)} = 1} I_n^{\cos} = 0. 
\end{equation*}
Here we have used $|{\psi} |_{H^1(\Omega)} = |\widetilde{\psi} |_{H^1(\Omega)} / c $. 
{This completes the argument for {\bf b)} and finishes the proof of Lemma \ref{lemma:constant}.}

\section{Concluding remarks}

This paper has originated as a complement to a part of \cite{Nie24} that deals with the preconditioned operator 
$\Delta^{-1}[\nabla \cdot (K\nabla u)]$ with the tensor $K(x, y)$ being diagonal and constant in an open subdomain $S\subset \Omega$. 
{However, we believe that present text goes beyond that: In addition to correcting the proof of Lemma 3.2 in the SIGEST paper \cite{Nie24} (and also in the original paper \cite{Ger20}), the new construction of the proof reveals an interesting relationship between the following two issues:}

\smallskip
\begin{description} 
    \item[I)] The preconditioned operator eigenvalue problem (i.e., the generalized eigenvalue problem (\ref{eq:eigen})) restricted to rectangular {\underline{subdomains} $\Sigma_l$ and constructed for a fixed eigenvalue $\lambda$. Here, $l$ denotes a parameter characterizing
    the size of $\Sigma_l$ that reduces to a single point in $S$ as $l \rightarrow 0$.}
    \item[II)] The spectrum of the preconditioned operator (the generalized eigenvalue problem (\ref{eq:eigen})) defined throughout the {\underline{entire}} domain $\Omega$.
    
\end{description}
\smallskip

Using the relationship with the wave equation, {\bf I)} provides for {a given $\lambda$, and an arbitrary small 
$\Sigma_l$,} the countably infinite set of eigenfunctions. When $\Sigma_l$ reduces to a single point in $S$, {this} $\lambda$ becomes a point in the spectrum {\bf II)}. This is how we should read the imprecise related comments on p.~133 and in the footnote on p.143 of \cite{Nie24}. 

The main result in \cite{Nie24} states that, under the assumption of symmetry and continuity of the tensor $K(x, y)$ throughout the closure $\overline{\Omega}$, the preconditioned operator spectrum {\bf II)} equals the convex hull of the ranges of diagonal functions in $\Lambda (x, y)$, {cf. the spectral decomposition presented in the Introduction}. Since this spectrum is an interval and 
the preconditioned operator is self-adjoint, this interval cannot be composed of eigenvalues. We therefore observe a transition of eigenvalues {of infinite multiplicity, associated with subdomains, to points of the continuous spectrum of the operator defined throughout the whole domain}. It remains unclear whether the spectral interval also contains any eigenvalues of the preconditioned operator defined throughout $\Omega$.   

The setting in this paper, as well as in~\cite{Ger20,Ger22,Nie24}, reflects the assumptions needed for the numerical solution of PDE boundary value problems. Therefore, the operators considered, {acting on Sobolev spaces,} are bounded and coercive. Consequently, the preconditioned operator is not compact and the associated generalized eigenvalue problem cannot be correctly defined {nor its solution numerically approximated} without further analysis justifying the {\em existence of eigenvalues}. In fact, the construction of an approximating sequence $(u_n)$ in Section \ref{sec:main_result} heavily relies on the non-compactness of the preconditioned operator. On the one hand, the sequence converges weakly to zero, while on the other hand, its norm remains bounded below by a positive constant.
One should notice that the mathematical physics and operator theory literature can use different settings and assumptions, which can lead to results that are not directly applicable here.

\end{document}